\documentclass[graybox]{svmult}

\usepackage{amsmath,amssymb,amsfonts}

\usepackage{mathptmx}

\begin{document}
	
\title*{On solution manifolds of some differential equations with more general state-dependent delay} 

\titlerunning{Solution manifolds for more general delay}
			
\author{Hans-Otto Walther}
	
\institute{Mathematisches Institut, Universit\"at Gie{\ss}en, Arndtstr. 2, Gie{\ss}en, 35392, Germany, \email{Hans-Otto.Walther@math.uni-giessen.de}}

\maketitle

\abstract {Differential equations with state-dependent delays define a semiflow of continuously differentiable solution operators in general only on the associated {\it solution manifold} in the Banach space $C^1_n=C^1([-h,0],\mathbb{R}^n)$. 
For a prototypic example we develop a new proof that its solution manifold is diffeomorphic to an open subset of the subspace given by $\phi'(0)=0$, without recourse to a restrictive hypothesis about the form of delays which is instrumental in earlier work on the nature of solution manifolds.
The new proof uses the framework of algebraic-delay systems.}

\medskip
	
	
\medskip


\begin{acknowledgement}
The author is grateful for a MATRIX workshop during which the approach used in this paper evolved close to its present state.
\end{acknowledgement}
 
\section{Introduction}

Let $h>0$. For $n\in\mathbb{N}$ let
$C_n=C([-h,0],\mathbb{R}^n)$ and $C^1_n=C^1([-h,0],\mathbb{R}^n)$ denote the Banach spaces of continuous and continuously dfferentiable maps $[-h,0]\to\mathbb{R}^n$, with the norms given by $|\phi|_{C_n}=\max_{-h\le t\le0}|\phi(t)|$ and $|\phi|_{C^1_n}=|\phi|_{C_n}+|\phi'|_{C_n}$, respectively. Abbreviate
$C=C_1$ and $C^1=C^1_1$. For maps with an interval $[t-h,t]$, $t\in\mathbb{R}$, in their domain of definition define $x_t$ by $x_t(s)=x(t+s)$ for $-h\le s\le0$.

\medskip

Differential equations with state-dependent delay define semiflows of continuously differentiable solution operators in general only on their solution manifolds. When written in the general form
$$
x'(t)=F(x_t)
$$
with a map $F:C^1_n\supset U_F\to\mathbb{R}^n$ then the solution manifold is defined as the set
$$
X_F=\{\phi\in U_F:\phi'(0)=F(\phi)\}.
$$
In \cite{W1,HKWW} it is shown that $X_F$ is a continuously differentiable submanifold of codimension $n$ in $C^1_n$ provided it is non-empty and $F$ is continuously differentiable with the additional smoothness property that

\medskip

(e) {\it each derivative $DF(\phi):C^1_n\to\mathbb{R}^n$, $\phi\in U_F$, has a linear extension $D_eF(\phi):C_n\to\mathbb{R}^n$ so that the map}
$$
U_F\times C_n\ni(\phi,\chi)\mapsto D_eF(\phi)\chi\in\mathbb{R}^n
$$
{\it is continuous.}

\medskip

See \cite{M-PNP} for a first version of property (e).

\medskip

What can be said about the submanifold $X_F$ ? The proof of \cite[Lemma 1]{KR} shows how to write $X_F$ as a graph with respect to a decomposition 
$$
C^1_n=X_0\oplus Q,\quad\dim\,Q=n,
$$
with the trivial solution manifold
$$
X_0=\{\phi\in C^1_n:\phi'(0)=0\},
$$
under the hypothesis that the extended derivatives $D_eF(\phi)\in L_c(C_n,\mathbb{R}^n)$ are bounded. A result in \cite{W6} guarantees a graph representation as before if delays are bounded away from zero, which in terms of $F$ means that the values $F(\phi)$ do not depend on $\phi(t)$ for $t$ in some fixed interval $[-h_F,0]$ with $-h<-h_F<0$. An example in \cite[Section 3]{W6} shows that without delays bounded away from zero solution manifolds do in general not
admit any graph representation in $C^1_n$ whatsoever. However, in \cite{W6} it is shown that for a large class of systems 
$$
x'(t)=g(x(t-d_1(x_t),\ldots,x(t-d_k(X_t)))\in\mathbb{R}^n
$$
with negative delays solution manifolds are {\it almost graphs}, which implies that they are diffeomorphic to open subsets of $X_0$. A crucial hypothesis in \cite{W6} is that the delays $d_{\kappa}$ admit a factorization
$$
d_{\kappa}(\phi)=\delta_{\kappa}(L\phi)
$$
with real functions $\delta_{\kappa}$ from an open subset $W$ of a finite-dimensional vectorspace $V$ into the interval $(0,h]$, and with a continuous linear map $L:C_n\to V$. In \cite{KW2} a proof which differs from its precursor in \cite{W6} in an essential part yields almost graph representations for systems as considered in \cite{W6} but without the restriction $\delta_{\kappa}(w)>0$.

\medskip

The recent paper \cite{W10} extends the result of \cite{KW2} to {\it algebraic-delay systems} \cite{W4} which consist of a differential equation with delays $r_{\kappa}(t)$ together with an algebraic equation which relates the  state $x_t\in C^1_n$ and the delay vector $r(t)\in\mathbb{R}^k$ to each other.
The hypotheses used in \cite{W10} include a factorization of the algebraic part of the system which involves linear maps into a finite-dimensional vectorspace,

\medskip

The said factorization hypotheses are not satisfied in various examples which arise from applications or are of interest for understanding better the impact of delays on dynamics in general, see e. g. the delay differential equations studied in \cite{GHMWW,W2,W3,L-WW,W5}.

\medskip

In the present case study we avoid a factorization hypothesis on the delay. We consider the prototypic equation 
\begin{equation}
x'(t)=f(x(t+d(x_t)))
\end{equation}
for $f:\mathbb{R}\to\mathbb{R}$ and 
$d:C^1\to(-h,0)$ both continuously differentiable, with the extension property that 

\medskip

(e,d) {\it each derivative $Dd(\phi):C^1\to\mathbb{R}$, $\phi\in C^1$, has a linear extension $D_ed(\phi):C\to\mathbb{R}$ so that the map
$$
C^1\times C\ni(\phi,\chi)\mapsto D_ed(\phi)\chi\in\mathbb{R}
$$
is continuous}.

\medskip

Eq. (1) takes the form $x'(t)=F(x_t)$ for $F:C^1\to\mathbb{R}$ given by $F(\phi)=f(\phi(d(\phi)))=f(ev(\phi,d(\phi)))$, with the evaluation map 
$$
ev:C^1\times(-h,0)\ni(\phi,t)\mapsto\phi(t)\in\mathbb{R}
$$
which is continuously differentiable with
$$
D\,ev(\phi,t)(\hat{\phi},\hat{t})=\hat{\phi}(t)+\phi'(t)\hat{t}.
$$

\begin{proposition} 
If condition (e,d) holds then the map $F$ is continuously differentiable and has property (e), and $X_F\neq\emptyset$.
\end{proposition}

\begin{proof} 1. Continuous differentiability of $F$ follows by means of the chain rule. For $\phi$ and $\chi$ in $C^1$
we get
$$
DF(\phi)\chi=f'(\phi(d(\phi)))[\chi(d(\phi))+\phi'(d(\phi))Dd(\phi)\chi].
$$
Define $D_eF(\phi):C\to\mathbb{R}$ by the previous formula with $Dd(\phi)$ replaced by $D_ed(\phi)$. Use the continuity of 
the evaluation map $ev_C:C\times(-h,0)\to\mathbb{R}$, $ev_C(\chi,r)=\chi(r)$, in combination with the continuity of differentiation $C^1\ni\phi\mapsto\phi'\in C$ in order to verify property (e). 

2. Proof of $X_F\neq\emptyset$. The set 
$K=\{\phi\in C^1:|\phi|_C\le 1\}$ is convex, and
$$
c=\sup_{\phi\in K}|f(\phi(-d(\phi)))|\le\sup_{|\xi|\le1}|f(\xi)|<\infty.
$$
Choose  $\phi_{\pm}\in K$ with 
$$
\phi_-'(0)<-c,\quad\phi_+'(0)>c.
$$
The map
$$
\Phi:[0,1]\ni t\mapsto \phi_-+t(\phi_+-\phi_-)\in K\subset C^1
$$
is continuous, hence the function
$$
z:[0,1]\ni t\mapsto\Phi(t)'(0)-f(\Phi(t)(d(\Phi(t))))\in\mathbb{R}
$$
is continuous as well. By the choice of $\phi_{\pm}$,
$z(0)<0<z(1)$, hence there is a zero $s$ of $z$. Then $\phi=\Phi(s)$
satisfies $\phi'(0)-F(\phi)=0$ which means $\phi\in X_F$; $X_F\neq\emptyset$. \qed
\end{proof}

It follows that $X_F$ is a continuously differentiable submanifold of $C^1$, with codimension 1.

\medskip

In the sequel we deal with the restrictions of $F$ to the open sets
$$
U_b=\{\phi\in C^1:|\phi'(t)|<b\,\,\mbox{on}\,\,[-h,d(\phi)]\}, \quad b>0.
$$
Notice that $C^1=\cup_{b>0}U_b$. Each non-empty solution manifold
$$
X_{F,b}=\{\phi\in U_b:\phi'(0)=F(\phi)\}=X_F\cap U_b
$$
associated with the restriction of $F$ to $U_b$ is a continuously differentiable submanifold of $C^1$ with codimension 1. Now we can state the main result of the present paper.

\begin{theorem}
Let $b>0$ with $X_{F,b}\neq\emptyset$ be given and assume
$$
sup_{\{\phi\in C^1:|\phi'(-h)|<b\}}|D_ed(\phi)|_{L_c(C,\mathbb{R})}<\infty.
$$	
Then there exists a diffeomorphism from $U_b$ onto an open subset of $C^1$ which maps the solution manifold $X_{F,b}$ onto an open subset of $X_0$.	
\end{theorem}

The hypothesis in Theorem 1.2 is a bit weaker that boundedness of $|D_ed(\phi)|_{L_c(C,\mathbb{R})}$ on all of $C^1$.

\medskip

The proof of Theorem 1 is inspired by the extension of work in \cite{W6,KW2} to algebraic-delay systems in \cite{W10} but uses an approach in its key part which is different from its counterparts in \cite{W6,KW2,W10}. We first embed $C^1$ as the graph of $d$ into the space $C^1\times\mathbb{R}$, by the graph map $\phi\stackrel{\Gamma}{\mapsto}(\phi,d(\phi))$.
For every $b>0$ let
$$
U_b^{\times}=\{(\phi,r)\in C^1\times(-h,0):|\phi'(t)|<b\,\,\mbox{on}\,\,[-h,r]\}.
$$
The sets $U_b^{\times}$, $b>0$, are open with 
$\Gamma(U_b)\subset U_b^{\times}$.

\medskip

Incidentally, we mention that $\Gamma(X_{F,b})$, $b>0$, is the solution manifold of the algebraic-delay system
\begin{align}
\phi'(0) & = f(\phi(r)),\nonumber\\
0 & = d(\phi)-r.\nonumber
\end{align}
in the set $U_b^{\times}$.

\medskip

In Section 2 below we find diffeomorphisms $T:C^1\times(-h,0)
\to C^1\times(-h,0)$ with $T(\Gamma(X_F))\subset X_0\times(-h,0)$ and $T(U_b^{\times})=U_b^{\times}$ for every $b>0$. The diffeomorphisms $T$ rely on continuously differentiable families of transversals to the space $X_0$. Section 3 contains the construction of such a family as it is needed for the proof of Theorem 1, which follows in Section 4.
The core of the proof is to show that for every $b>0$ the projection $P:C^1\times\mathbb{R}\ni(\phi,r)\mapsto\phi\in C^1$ defines a diffeomorphism from the submanifold $T(\Gamma(U_b))$ of codimension 1 in $C^1\times\mathbb{R}$ onto an open subset of $C^1$. This, in turn, is achieved by an estimate of the partial derivative $D_2$ of the first component of the diffeomorphism $T$ given by the transversals constructed in Section 3.

\medskip

Section 5
provides examples which show that the hypotheses of Theorem 1 do not imply any of the sufficient conditions from \cite{KR,W6,KW2} for a graph or almost graph representation of the solution manifolds $X_{F,b}\neq\emptyset$.

\medskip

It is an open problem whether the detour towards the desired diffeomorphism via the solution manifold of an algebraic-delay system can be avoided. 
Among further open problems we mention removing the hypothesis $d(\phi)<0$, and generalization to larger systems with several discrete delays. For the latter, as well as for allowing $d(\phi)=-h$, the papers \cite{W6,KW2,W10} provide tools and techniques. 

\medskip

{\bf Notation, preliminaries.} The closure of a subset $S$ of a topological space $T$ is denoted by $cl\,S$. The relation $A\subset\subset B$ for open subsets $A,B$ of $T$ means that $cl\,A$ is compact and contained in $B$.

\medskip 

For subsets $A,B$ of a vectorspace $V$ over a field $\mathbb{K}$ and for 
$x\in\mathbb{K},\,\,M\subset\mathbb{K}$ the sets $A+B$, $aB$, $MB$ are defined in the obvious way, e. g.,
$$
A+B=\{c\in V:\mbox{For some}\,\,a\in A\,\,\mbox{and}\,\,b\in B,\,c=a+b\}.
$$

\medskip

Finite-dimensional vectorspaces are always equipped with the norm topology.

\medskip

The constant function on $[-h,0]$ with value $y$ is denoted by $\mathbf{y}$.

\medskip

Derivatives and partial derivatives are continuous linear maps and indicated by capitals $D$ and $D_j$, respectively. For maps on domains $dom\subset\mathbb{R}$, $\phi'(t)=D\phi(t)1$.

\section{A diffeomorphism $C^1\times(-h,0)\to C^1\times(-h,0)$}

There exist continuously differentiable maps 
$$
\chi:\mathbb{R}\times(-h,0)\to C^1
$$
which satisfy
\begin{align}
\chi(v,r)'(0) & = 1\\
& \mbox{and} \nonumber\\
\chi(v,r)(r) & = 0\quad\mbox{for all}\quad(v,r)\in\mathbb{R}\times(-h,0),
\end{align}
for example, the affine linear map given by  $\chi(v,r)(t)=t-r$.
By (2), $\chi(v,r)\notin X_0$ and $C^1=X_0\oplus\mathbb{R}\,\chi(v,r)$ for every $(v,r)\in\mathbb{R}\times(-h,0)$ (transversality).

\begin{proposition}
Suppose $\chi:\mathbb{R}\times(-h,0)\to C^1$ is continuously differentiable and satisfies (2). Then the map $A:C^1\times(-h,0)\to C^1$ given by 
$$
A(\phi,r)=\phi-f(\phi(r))\,\chi(\phi(r),r)
$$
is continuously differentiable, we have
$$
A(\Gamma(X_F))\subset X_0,
$$
and the map 
$$
T:C^1\times(-h,0)\ni(\phi,r)\mapsto(A(\phi,r),r)\in C^1\times(-h,0)
$$ 
is a continuously differentiable diffeomorphism onto $C^1\times(-h,0)$. The inverse $Y:C^1\times(-h,0)\to C^1\times(-h,0)$ of $T$ is given by
$$
Y(\psi,r)=(B(\psi,r),r)
$$ 
with the continuously differentiable map 
$$
B:C^1\times(-h,0)\ni(\psi,r)\mapsto\psi+f(\psi(r))\,\chi(\psi(r),r)\in  C^1.
$$ 
\end{proposition}

\begin{proof} 1. The maps $A,T,B,Y$ are continuously differentiable.

\medskip

2. Let $\phi\in X_F$ be given. Using (2) we get  
\begin{align}
A(\Gamma(\phi))'(0) & = A(\phi,d(\phi))'(0)=\phi'(0)-f(\phi(d(\phi)))\,\chi(\phi(d(\phi)),d(\phi))'(0)\nonumber\\
& =  \phi'(0)-f(\phi(d(\phi)))=0,\nonumber
\end{align}
which proves $A(\Gamma(\phi))\in X_0$.

\medskip

3. Proof of $Y(T(\phi,r))=(\phi,r)$ everywhere: For $(\phi,r)\in C^1\times(-h,0)$ set $\psi=A(\phi,r)$. Then
$$
Y(T(\phi,r))=Y(A(\phi,r),r)=Y(\psi,r)=(B(\psi,r),r)
$$
and
$$
B(\psi,r)=\psi+f(\psi(r))\,\chi(\psi(r),r)=[\phi-f(\phi(r))\,\chi(\phi(r),r)]+f(\psi(r))\,\chi(\psi(r),r).
$$
From (3) we have
$$
\psi(r)=\phi(r)-f(\phi(r))\,\chi(\phi(r),r)(r)=\phi(r).
$$
It follows that $B(\psi,r)=\phi$, hence $Y(T(\phi,r))=(B(\psi,r),r)=(\phi,r)$.

\medskip

4. As in Part 3 one finds $T(Y(\psi,r))=(\psi,r)$ on $C^1\times(-h,0)$. Using this in combination with the  result of Part  3 we see that $T$ is a diffeomorphism onto $C^1\times(-h,0)$ whose inverse is the map $Y$. 
\qed
\end{proof}

\section{Transversals to $X_0$ which are small in $C$}

We need a map $\chi:\mathbb{R}\times(-h,0)\to C^1$ as in Section 2 which in addition to (2) and (3) satisfies further conditions. Let $I:C^1\hookrightarrow C$ denote the continuous linear inclusion map.

\begin{proposition}
Let a continuous function $H:\mathbb{R}\to(0,\infty)$ and $\epsilon>0$ be given. Then there exist continuously differentiable maps $\chi:\mathbb{R}\times(-h,0)\to C^1$ such that for every $(v,r)\in\mathbb{R}\times(-h,0)$ Eq. (2) holds, and
\begin{align}
\chi(v,r)(t) & =  0\quad\mbox{for all}\quad t\in[-h,r],\\
|\chi(v,r)|_C & \le  H(v),\\
|D(I\circ \chi)(v,r)|_{L_c(\mathbb{R}^2,C)} & \le  H(v).
\end{align}
\end{proposition}

\begin{proof} 1. For every $\epsilon>0$ and $z\in(-h,0)$ there exists $\psi\in C^1$ with
$$
\psi'(0)=1,\quad\psi(t)=0\,\,\mbox{for all}\,\,t\in[-h,z],\quad|\psi|_C<\epsilon.
$$

2. For $n\in\mathbb{N}$ set $U_n=(-n,n)\times\left(-h+\frac{h}{3n},-\frac{h}{3n}\right)$ and choose open sets $U_{n,0}\subset U_{n,1}$ so that for every $n\in\mathbb{N}$,
$$
U_n\subset\subset U_{n,0}\subset\subset U_{n,1}\subset\subset U_{n+1}.
$$
There exist continuously differentiable functions $a_n:\mathbb{R}\times(-h,0)\to[0,1]\subset\mathbb{R}$, $n\in\mathbb{N}$, so that for every $n\in\mathbb{N}$,
$$
a_n(v,r)=1\,\,\mbox{on}\,\,U_{n,0},\quad a_n(v,r)=0\,\,\mbox{on}\,\,\mathbb{R}\times(-h,0)\setminus cl\,U_{n,1}.
$$ 
Choose an increasing sequence $(A_n)_1^{\infty}$ in $(0,\infty)$ with
$$
A_n\ge\max_{v\in\mathbb{R},-h<r<0}|Da_n(v,r)|_{L_c(\mathbb{R}^2,\mathbb{R})}\,\,\mbox{for all}\,\,n\in\mathbb{N}.
$$
and set $H_n=\min_{-n-2\le v\le n+2}H(v)$, $n\in\mathbb{N}$. The sequence $(H_n)_1^{\infty}$ is decreasing, and
$$
\epsilon_n=\frac{H_n}{1+2A_n}
$$
defines a decreasing sequence of positive reals. For every $n\in\mathbb{N}$, Part 1 with $z=z_n=-\frac{h}{3(n+2)}$ and $\epsilon=\epsilon_n$ yields $\psi_n\in C^1$ with
$$
\psi_n'(0)=1,\quad\psi_n(t)=0\,\,\mbox{for all}\,\,t\in\left[-h,-\frac{h}{3(n+2)}\right],\quad|\psi_n|_C<\epsilon_n.
$$

3. We define continuously differentiable maps $\sigma_n:\mathbb{R}\times(-h,0)\supset dom_n\to C^1$, $n\in\mathbb{N}$, by
$$
dom_1=U_{2,0}\quad\mbox{and}\quad\sigma_1(v,r)=\psi_2,
$$
and for integers $n\ge2$,
$$
dom_n=U_{n+1,0}\setminus cl\,U_n\quad\mbox{and}\quad\sigma_n(v,r)=a_n(v,r)\psi_n+(1-a_n(v,r))\psi_{n+1}.
$$
The domains $dom_n$ overlap with
$$
dom_2\cap dom_1=(U_{3,0}\setminus cl\, U_2)\cap U_{2,0}=U_{2,0}\setminus cl\, U_2
$$
and
$$
dom_{n+1}\cap dom_n=(U_{n+2,0}\setminus cl\, U_{n+1})\cap( U_{n+1,0}\setminus cl\, U_n)=U_{n+1,0}\setminus cl\, U_{n+1}
$$
for integers $n\ge2$, while for integers $n\ge1$ and $j\ge n+2$ we have 
$$
dom_j\cap dom_n=\emptyset,
$$ 
due to
\begin{align}
dom_j\cap dom_1 & =  (U_{j+1,0}\setminus cl\, U_j)\cap U_{2,0}\nonumber\\
& =  U_{2,0}\setminus cl\,U_j\subset U_{3,0}\setminus cl\, U_3=\emptyset\nonumber
\end{align}
for $j\ge 3$ and
\begin{align}
dom_j\cap\, dom_n & =  (U_{j+1,0}\setminus cl\, U_j)\cap( U_{n+1,0}\setminus cl\, U_n)\nonumber\\
& =  U_{n+1,0}\setminus cl\, U_j\subset U_{n+1,0}\setminus cl\, U_{n+2}\subset U_{n+2}\setminus cl U_{n+2}=\emptyset\nonumber
\end{align}
for $n\ge2$ and $j\ge n+2$.

\medskip

4. Proof that $\sigma_n$ and $\sigma_{n+1}$ coincide on the overlap of their domains. On $dom_2\cap dom_1=U_{2,0}\setminus cl\,U_2$ we have $a_2(v,r)=1$ and thereby
$$
\sigma_2(v,r)=a_2(v,r)\psi_2+(1-a_2(v,r))\psi_3=\psi_2=\sigma_1(v,r).
$$
For integers $n\ge 2$ we have
$$
dom_{n+1}\cap dom_n=U_{n+1,0}\setminus cl\, U_{n+1},
$$
with $a_{n+1}(v,r)=1$ and $a_n(v,r)=0$ for $(v,r)\in U_{n+1,0}\setminus cl\, U_{n+1}$. For such $(v,r)$ we infer 
$$
\sigma_{n+1}(v,r)=a_{n+1}(v,r)\psi_{n+1}+(1-a_{n+1}(v,r))\psi_{n+2}=\psi_{n+1}
$$
and
$$
\sigma_n(v,r)=a_n(v,r)\psi_n+(1-a_n(v,r))\psi_{n+1}=\psi_{n+1},
$$
which means $\sigma_{n+1}(v,r)=\sigma_n(v,r)$.

\medskip

5. It follows that the maps $\sigma_n$ define a continuously differentiable map $\chi:\mathbb{R}\times(-h,0)\to C^1$.
In order to verify the conditions (2) and (4)-(6) let $(v,r)\in\mathbb{R}\times(-h,0)$ be given. Then $(v,r)\in dom_n$ for some integer $n\ge1$, and $\chi(v,r)=\sigma_n(v,r)$. From $dom_n\subset U_{n+1,0}\subset U_{n+2}$ we have $-n-2<v<n+2$ and $-h<r<-\frac{h}{3(n+2)}$.

\medskip

5.1 The case $n\ge 2$. Then
$$
\chi(v,r)'(0)=a_n(v,r)\psi_n'(0)+(1-a_n(v,r))\psi_{n+1}'(0)=1.
$$
For $-h\le t\le r<-\frac{h}{3(n+2)}<-\frac{h}{3((n+1)+2)}$, $\psi_n(t)=0$ and $\psi_{n+1}(t)=0$, hence
$$
\chi(v,r)(t)=a_n(v,r)\psi_n(t)+(1-a_n(v,r))\psi_{n+1}(t)=0.
$$
Next,
\begin{align}
|\chi(v,r)|_C & =  |a_n(v,r)\psi_n+(1-a_n(v,r))\psi_{n+1}|_C\nonumber\\
& \le  a_n(v,r)|\psi_n|_C+(1-a_n(v,r))|\psi_{n+1}|_C\nonumber\\
& \le  a_n(v,r)\epsilon_n+(1-a_n(v,r))\epsilon_{n+1}\le a_n(v,r)\epsilon_n+(1-a_n(v,r))\epsilon_n\nonumber\\
& =  \epsilon_n\le  H_n\le H(v)\nonumber
\end{align}
where the last inequality follows from $-n-2<v<n+2$. 

\medskip

The map $I\circ\chi$ is given on $dom_n$ by 
$$
(I\circ\sigma_n)(w,s)=a_n(w,s)I\psi_n+(1-a_n(w,s))I\psi_{n+1}.
$$ 
Using a product rule we differentiate and obtain the estimates 
\begin{align}
|D(I\circ\chi)(v,r)|_{L_c(\mathbb{R}^2,C)} & \le  |Da_n(v,r)|_{L_c(\mathbb{R}^2,\mathbb{R})}|\psi_n|_C+|Da_n(v,r)|_{L_c(\mathbb{R}^2,\mathbb{R})}|\psi_{n+1}|_C\nonumber\\
& \le  A_n\epsilon_n+A_n\epsilon_{n+1}\le2A_n\epsilon_n\le H_n\le H(v).\nonumber
\end{align}

\medskip

5.2 The case $n=1$. As $\sigma_1$ is constant on $dom_1=U_{2,0}\subset U_3$ with value $\psi_2$ we infer $\chi(v,r)'(0)=\psi_2'(0)=1$. From $-3<v<3$ we get $|\chi(v,r)|_C=|\psi_2|_C<\epsilon_2<H_2\le H(v)$. Obviously, $D(I\circ\chi)(v,r)=0$.
For $-h\le t\le r<-\frac{h}{3\cdot 3}<-\frac{h}{3(2+2)}$, $\chi(v,r)(t)=\psi_2(t)=0$. \qed
\end{proof}

\section{Proof of Theorem 1}

1. Let $b>0$ be given and set
$$
c=sup_{\{\phi\in C^1:|\phi'(-h)|<b\}}|D_ed(\phi)|_{L_c(C,\mathbb{R})}<\infty.
$$
Let $H:\mathbb{R}\to(0,\infty)$ denote the continuous function given by 
$$
H(v)=\frac{1}{2(b+1)c(1+|f(v)|+|f'(v)|)}
$$ 
and apply Proposition 3. This results in a map $\chi$ for which in particular conditions (2) and (3) hold. Therefore Proposition 2 applies. Consider the continuously differentiable maps $A,T,B,Y$ according to Proposition 2, which now involve $\chi$ obtained from Proposition 3, with properties (2) and (4)-(6). 

\medskip

2. Proof of $T(U_b^{\times})=U_b^{\times}$. 

\medskip

2.1 Let $(\phi,r)\in U_b^{\times}$. For all $t\in[-h,r]$,
$$
A(\phi,r)(t)=\phi(t)-f(\phi(r))\chi(\phi(r),r)(t)=\phi(t)
$$
since $\chi(\phi(r),r)(t)=0$ due to property (4). Using $(\phi,r)\in U_b^{\times}$ we obtain $|A(\phi,r)'(t)|=|\phi'(t)|<b$ for all $t\in[-h,r]$, which means
$T(\phi,r)=(A(\phi,r),r)\in U_b^{\times}$. It follows that 
$T(U_b^{\times})\subset U_b^{\times}$.

\medskip

2.2 Analogously one gets $Y(U_b^{\times})\subset U_b^{\times}$, which yields $U_b^{\times}\subset T(U_b^{\times})$. 

\medskip

3. Proof of
$$
|D_2(d\circ B)(\psi,r)1|\le\frac{1}{2}\quad\mbox{for every}\quad(\psi,r)\in U_b^{\times}.
$$
At $(\psi,r)\in U_b^{\times}$,
\begin{align}
D_2(d\circ B)(\psi,r)1 & =  Dd(B(\psi,r))DB(\psi,r)(0,1)\nonumber\\
& =  D_ed(B(\psi,r)) \,I\,DB(\psi,r)(0,1)\nonumber\\
& =  D_ed(B(\psi,r))D(I\circ B)(\psi,r)(0,1)\nonumber\\
& = D_ed(B(\psi,r))D_2(I\circ B)(\psi,r)1\nonumber
\end{align}
with
$$
D_2(I\circ B)(\psi,r)1=f'(\psi(r))\psi'(r)\,I\,\chi(\psi(r),r)+f(\psi(r))\,D(I\circ\chi)(\psi(r),r)(\psi'(r),1).
$$
By Part 2, $Y(\psi,r)\in U_b^{\times}$, since  $(\psi,r)\in U_b^{\times}$. It follows that $|B(\psi,r)'(-h)|<b$. Using the bound $c$ for $|D_ed(\phi)|_{L_c(C,\mathbb{R})}$ on $\{\phi\in C^1:|\phi'(-h)|<b\}$, and $|\psi'(r)|<b$ due to $(\psi,r)\in U_b^{\times}$, we get
\begin{align}
|D_2(d\circ B)(\psi,r)1| & \le  |D_ed(B(\psi,r))|_{L_c(C,\mathbb{R})}\{|f'(\psi(r))||\psi'(r)||\,I\,\chi(\psi(r),r)|_C\nonumber\\
&  +|f(\psi(r))||D(I\circ\chi)(\psi(r),r)|_{L_c(\mathbb{R}^2,C)}|(\psi'(r),1)|\}\nonumber\\
& \le  c\{|f'(\psi(r))|\,b\,H(\psi(r))+|f(\psi(r))|\,H(\psi(r))|(b+1)\}\nonumber\\
& \le  c(b+1)\,H(\psi(r))(1+|f'(\psi(r))|+|f(\psi(r))|)=\frac{1}{2}.\nonumber
\end{align}

\medskip

4. From $T(\Gamma(U_b))\subset T(U_b^{\times})=U_b^{\times}$, we infer
\begin{align}
T(\Gamma(U_b)) & = \{(\psi,r)\in U_b^{\times}:Y(\psi,r)\in\Gamma(U_b))\}\nonumber\\
& = \{(\psi,r)\in U_b^{\times}:B(\psi,r)\in U_b\,\,\mbox{and}\,\,0=d(B(\psi,r))-r\}.\nonumber
\end{align}
Notice that $\{(\psi,r)\in U_b^{\times}:B(\psi,r)\in U_b\}\subset C^1\times(-h,0)$ is open in $C^1\times\mathbb{R}$.

\medskip

5. Using Part 3 we apply the Implicit Function Theorem at the points $(\psi,r)\in T(\Gamma(U_b))$ to the equation $0=d(B(\psi,r))-r$  and obtain that locally the zeroset $T(\Gamma(U_b))$  is given by continuously differentiable functions from open subsets of $C^1$ into $(-h,0)$. This implies that the image of $T(\Gamma(U_b))$ under the projection $P:C^1\times\mathbb{R}\ni(\psi,r)\mapsto\psi\in C^1$ is an open subset of $C^1$.

\medskip

6. Proof that $P$ is injective on $T(\Gamma(U_b))$. Let $(\psi,r)$ and $(\phi,s)$ in $T(\Gamma(U_b))\subset U_b^{\times}$ be given with $P(\psi,r)=P(\phi,s)$. Then $\psi=\phi$. It remains to show $r=s$. 

\medskip

6.1 We show that the line segment $(\psi,r)+[0,1]((\psi,s)-(\psi,r))=(\psi,s)+[0,1]((\psi,r)-(\psi,s))$ is contained in $U_b^{\times}$.

\medskip

6.1.1 In case $r\le s$ we have $(\psi,r)+[0,1]((\psi,s)-(\psi,r))=\{\psi\}\times[r,s]$.
Let $t\in[r,s]$ be given. We need to show $(\psi,t)\in U_b^{\times}$. For each $\tau\in[-h,t]$,  $\tau\in[-h,s]$. Using $(\psi,s)=(\phi,s)\in U_b^{\times}$ we get
$$
|\psi'(\tau)|=|\phi'(\tau)|<b.
$$ 
It follows that $(\psi,t)\in U_b^{\times}$. 

\medskip

6.1.2 In case $s\le r$ we have $(\psi,r)+[0,1]((\psi,s)-(\psi,r))=\{\psi\}\times[s,r]$.
Let $t\in[s,r]$ be given. We need to show $(\psi,t)\in U_b^{\times}$. For each $\tau\in[-h,t]$,  $\tau\in[-h,r]$. Using $(\psi,r)\in U_b^{\times}$ we have
$$
|\psi'(\tau)|<b.
$$ 
It follows that $(\psi,t)\in U_b^{\times}$. 

\medskip

6.2 Using the equation defining $T(\Gamma(U_b))$ and Part 6.1 we get
$$
s-r=d(B((\psi,s))-d(B((\psi,r))=\int_0^1D_2(d\circ B)(\psi,r+\theta(s-r))[s-r]d\theta.
$$
Applying the estimate in Part 3 we arrive at $|s-r|\le\frac{1}{2}|s-r|$, hence $s=r$. 

\medskip

7. The results of Parts 5 and 6 combined yield that the equation $\chi=PT(\Gamma(\phi))$ defines a diffeomorphism  from the open subset $U_b\subset C^1$ onto the open subset $V_b=PT(\Gamma(U_b))\subset C^1$. This diffeomorphism maps the submanifold $X_{F,b}\subset U_b$ onto a continuously differentiable submanifold $PT(\Gamma(X_{F,b}))$ of codimension 1 in $C^1$. From Proposition 2 we obtain that $PT(\Gamma(X_{F,b}))$  is contained in the closed subspace $X_0$ of codimension 1 in $C^1$. We infer that $PT(\Gamma(X_{F,b}))$ is an open subset of the space $X_0$. \qed

The details of the last argument in Part 7 of the proof are as follows, for a closed subspace $V$ of finite codimension in a Banach space $B$ and a continuously differentiable submanifold $X\subset V$ of $B$, with the same codimension in $B$: From $X\subset V$, $T_xX\subset V$ for all $x\in X$.
As both spaces have the same finite codimension, $T_xX=V$.
The restriction of $id:B\to B$ defines a continuously differentiable map $\xi:X\to V$ with $D\xi(x)z=z$ on $T_xV=V$ for each $x\in X$. Hence each derivative $D\xi(x)$ is an isomorphism. It follows that $X=\xi(X)$ is an open subset of $V$.  

\section{Examples}

We give examples of continuously differentiable maps $f:\mathbb{R}\to\mathbb{R}$ and $d:C^1\to(-h,0)$ which satisfy the hypotheses for Theorem 1 but none of the hypotheses in \cite{KR,W6,KW2} for graph or almost graph representations of the solution manifolds $X_{F,b}$  in the sets $U_b$, $b>0$. 

\medskip

We begin with violation of the hypothesis from \cite{KR} for graph representations of the solution manifolds $X_{F,b}$, $b>0$, which is boundedness of the extended derivatives $D_eF(\phi)$ on the sets $U_b$.

\begin{proposition}
In case condition (e,d) holds and
\begin{equation}
\sup_{n\in\mathbb{N}}|f'(n)|=\infty,
\end{equation}
$$
\sup_{\phi\in U_b}|D_eF(\phi)|_{L_c(C,\mathbb{R})}=\infty\quad\mbox{
for every}\quad b>0.
$$
\end{proposition}

\begin{proof} Let $b>0$ be given. For every $n\in\mathbb{N}$, obviously $\mathbf{n}\in U_b$.
We have
\begin{align}
	|D_eF(\mathbf{n})|_{L_c(C,\mathbb{R})} & \ge  |D_eF(\mathbf{n})\mathbf{1}|=|DF(\mathbf{n}){\mathbf 1}|\nonumber\\	& =  |Df(\mathbf{n}(d(\mathbf{n})))Dev(\mathbf{n},d(\mathbf{n}))({\mathbf 1},Dd(\mathbf{n}){\mathbf 1})|\nonumber\\
	& =  |f'(\mathbf{n}(d(\mathbf{n})))\{{\mathbf 1} (d(\mathbf{n}))+\mathbf{n}'(d(\mathbf{n}))\,Dd(\mathbf{n}){\mathbf 1}\}|\nonumber\\	& =  |f'(n)|\{1+0\}=|f'(n)|,\nonumber
\end{align}
hence
$$
\sup_{\phi\in U_b}|D_eF(\phi)|_{L_c(C,\mathbb{R})}\ge\sup_{n\in\mathbb{N}}|f'(n)|,
$$
which shows that $|D_eF(\phi)|_{L_c(C,\mathbb{R})}$ is unbounded on $U_b$ in case $\sup_{n\in\mathbb{N}}|f'(n)|=\infty$. 
\qed
\end{proof}

We turn to violation of the hypothesis from \cite[Theorem 2.4]{W6} for graph representations of the solution manifolds $X_{F,b}$. Given $b>0$ the condition of interest for the map $U_b\ni\phi\mapsto F(\phi)\in\mathbb{R}$  is that delays are bounded away from zero in the following sense: There exists $s_b\in(-h,0)$
so that $F(\phi)=F(\psi)$ for all $\phi,\psi$ in $U_b$
with $\phi(t)=\psi(t)$ for all $t\in[-h,s_b]$.

\begin{proposition}
(i) Assume $v:\mathbb{R}\to\mathbb{R}$ and $\delta:\mathbb{R}\to(-h,0)$ are continuously differentiable  with bounded derivatives. Then the map $d:C^1\to(-h,0)$ given by
\begin{equation}
d(\phi)=\delta\left(\int^0_{-h}v(\phi(t))dt\right)
\end{equation}
is continuously differentiable with property (e,d) and satisfies
\begin{equation}
\sup_{\phi\in C^1}|D_ed(\phi)|_{L_c(C,\mathbb{R})}\,\,<\,\,\infty.
\end{equation}
(ii) Assume $f:\mathbb{R}\to\mathbb{R}$ is continuously differentiable and injective, and $v$ and $\delta$ are given as in assertion (i), with the additional properties
\begin{equation}
\lim_{y\to\infty}v(y)=\infty\quad\mbox{and}\quad\lim_{w\to\infty}\delta(w)=0.
\end{equation}
Consider the delay functional
$d:C^1\to(-h,0)$ defined by Eq. (8), and $F:C^1\to\mathbb{R}$
defined by $F(\phi)=f(\phi(d(\phi)))$. Let $b>0$. For every $s\in(-h,0)$ there exist $\phi$ and $\psi$ in $U_b$ with $\phi(t)=\psi(t)$ for all $t\in[-h,s]$ and
$F(\phi)\neq F(\psi)$.
\end{proposition}

\begin{proof} 1. On (i). By \cite[Lemma 1.5, Appendix IV]{DvGVLW} the substitution operator $V:C\ni\phi\mapsto v\circ\phi\in C$ is continuously differentiable with $(DV(\phi)\chi)(t)=v'(\phi(t))\chi(t)$ for $\phi,\chi$ in $C$ and $t\in[-h,0]$. Integration $C\to\mathbb{R}$ is linear and continuous. The chain rule yields that the map $d_C:C\to(-h,0)$
given by $d_C(\phi)=\delta\left(\int^0_{-h}v(\phi(t))dt\right)$ 
is continuously differentiable. This implies that also the map $d:C^1\to(-h,0)$ is continuously differentiable, and it follows that condition $(e,d)$ is satisfied with $D_ed(\phi)=Dd_C(\phi)$ for all $\phi\in C^1$. 
Boundedness of $D_ed(\phi)$, $\phi\in C^1$, becomes obvious from
$$
|(Dd_C(\phi)\chi)|\le\left|\delta'\left(\int_{-h}^0V(\phi)\right)\right|
\left|\int_{-h}^0DV(\phi)\chi\right|\le\sup_{w\in\mathbb{R}}|\delta'(w)\left|\int_{-h}^0DV(\phi)\chi\right|
$$ 
with $|DV(\phi)\chi|_C\le\sup_{y\in\mathbb{R}}|v'(y)||\chi|_C$ for all $\phi\in C$ and all $\chi\in C$. 

\medskip

2. On (ii). Let $b>0$ and $s\in(-h,0)$ be given. Choose $w_s>0$ with $\delta(w)>s$ on $[w_s,\infty)$. Choose $c>0$ so large that $v(\tilde{c})>\frac{w_s}{h}$ for $\tilde{c}\ge c$. Consider $\phi={\bf c}$ and choose $\psi\in C^1$ with
$$
\psi(t)=c\,\,\mbox{on}\,\,[-h,s],\,\, c<\psi(t)\,\,\mbox{on}\,\,(s,0],\,\, |\psi'(t)|<b\,\,\mbox{on}\quad(s,0].
$$
Then both $\phi$ and $\psi$ belong to $U_b$, and $\phi(t)=\psi(t)$ on $[-h,s]$. Obviously, $F(\phi)=f(c)$.
It remains to show $F(\psi)\neq f(c)$. From $\psi(t)\ge c$ on $[-h,0]$  we have $v(\psi(t))\ge\frac{w_s}{h}$ for all $t\in[-h,0]$, hence
$$
\int_{-h}^0V\circ\psi\ge h\frac{w_s}{h}=w_s,
$$
which yields 
$$
d(\psi)=\delta\left(\int_{-h}^0V\circ\psi\right)>s.
$$
Consequently, $\psi(d(\psi))>c$. By the injectivity of $f$,
$$
F(\phi)=f(c)\neq f(\psi(d(\psi)))=F(\psi).\quad\qed
$$
\end{proof}

Finally we look for delay functionals $d$ which do not admit factorizations as they were assumed in the results \cite[Theorem 5.1]{W6} and \cite[Theorem 3.5]{KW2} on almost graph 
representations of solution manifolds. More exactly, on each set $U_b$, $b>0$, the desired delay functionals should not have the form
\begin{equation}
d(\phi)=\hat{\delta}(L\phi)
\end{equation}
with a continuous linear map $L:C\to V$ into a finite-dimensional vectorspace $V$ and a continuously differentiable function $\hat{\delta}:W\to(-h,0)$, $W\subset V$ open. The form (11) for $d$ on $U_b$ implies $U_b\subset L^{-1}(W)\cap C^1$ and $d(\phi)=\hat{\delta}(0)$ on $U_b\cap (L^{-1}(0)\cap C^1)$. 
The subspace $Z=L^{-1}(0)\cap C^1=(L|C^1)^{-1}(0)\subset C^1$ is closed and has finite codimension in $C^1$. We rephrase:
If $d$ has the form (11) then it is constant on $U_b\cap Z$ for some closed subspace $Z\subset C^1$ of finite codimension. The next proposition provides delay functionals for which the previous necessary condition is violated. 

\begin{proposition}
Assume $v:\mathbb{R}\to\mathbb{R}$ and $\delta:\mathbb{R}\to(-h,0)$ are continuously differentiable  with bounded derivatives. Assume in addition that $\delta$ is injective and
$$
v(y)=0\,\,\mbox{on}\,\,(-\infty,0],\,\,v'(y)>0\,\,\mbox{on}\,\,(0,\infty).
$$
Let $d:C^1\to(-h,0)$ be defined by Eq. (8), let $b>0$, and let a closed subspace $Z\subset C^1$ of finite codimension be given. Then $d$  is not constant on $U_b\cap Z$.
\end{proposition}

\begin{proof} Let a closed subspace $Z\subset C^1$ of finite codimension be given, and let $b>0$. Obviously, $Z\neq\{0\}$. Choose $\phi\in Z\setminus\{0\}$. Multiplying with a sufficiently small real number if necessary we achieve $|\phi'(t)|<\frac{b}{2}$ on $[-h,0]$, and $\phi(t)>0$ for some $t\in[-h,0]$. Both $\phi$ and $2\phi$ belong to $U_b\cap Z$. Using the properties of $v$ we get
\begin{align}
\int_{-h}^0v(\phi(t))dt & =  
\int_{\{t\in[-h,0]:\phi(t)>0\}}v(\phi(t))dt\nonumber\\
& <  \int_{\{t\in[-h,0]:\phi(t)>0\}}v(2\phi(t))dt=\int_{-h}^0v(2\phi(t))dt,\nonumber
\end{align}
and injectivity of $\delta$ yields
$$
d(\phi)=\delta\left(\int_{-h}^0v(\phi(t))dt\right)\neq
\delta\left(\int_{-h}^0v(2\phi(t))dt\right)=d(2\phi).\quad\qed
$$
\end{proof}

Notice that the hypotheses on $d:C^1\to\mathbb{R}$ and  $f:\mathbb{R}\to\mathbb{R}$ which are required in Propositions 4-6
are compatible. Consequently there exist continuously differentiable maps
$d:C^1\to(-h,0)$ and $f:\mathbb{R}\to\mathbb{R}$  which satisfy the hypotheses of Theorem 1 while 
none of the results from \cite{KR,W6,KW2}
on graph and almost graph representations applies to the
solution manifolds $X_{F,b}$ of the restrictions $F|U_b$, $b>0$.

\medskip

A caveat remains: The hypotheses which are violated by the previous examples concern the restrictions of $F$ to the open sets $U_b$, $b>0$, which are neighbourhoods of the solution manifolds $X_{F,b}$. For the previous examples we did not exclude the possibility that for some $b>0$ on some neighbourhood of $X_{F,b}$ which is smaller than $U_b$  results from \cite{KR,W6,KW2} can be applied to $X_{F,b}$.

\end{document}